\newcommand{\Posa}{P{\'o}sa}
\newcommand{\mlevels}{middle levels}

\documentclass{elsart}

\usepackage{amssymb}
\usepackage{amsmath}
\usepackage{latexsym}
\usepackage{graphicx}
\usepackage{subfigure}          

\newtheorem{theorem}{Theorem}

\begin{document}

\begin{frontmatter}



\title{An update on the \mlevels\ problem}


\author{Ian Shields\corauthref{cor}}
\corauth[cor]{Corresponding author.}
\address{IBM, P.O. Box 12195, 
Research Triangle Park, North Carolina 27709, USA}
\ead{ishields@us.ibm.com}

\author{Brendan J. Shields}
\address{Department of Physics, 
Massachusetts Institute of Technology, 
77 Massachusetts Avenue
Cambridge, MA 02139-4307, USA}
\ead{bshields@fas.harvard.edu}

\author{Carla D. Savage\thanksref{cdsthanks}}
\address{Department of Computer Science, 
North Carolina State University, Box 8206,
Raleigh, North Carolina 27695-8206, USA}
\ead{savage@cayley.csc.ncsu.edu}

\thanks[cdsthanks]{Research supported in part by NSF grant DMS-0300034}
\begin{abstract}
The \mlevels\ problem is to find a Hamilton cycle in the
\mlevels, $M_{2k+1}$, of the Hasse diagram of $ {\mathcal
B}_{2k+1} $ (the partially ordered set of subsets of a $2k+1$-element
set ordered by inclusion). Previously, the best known, from
\cite{MR2000k:05179}, was that $ {M}_{2k+1}$ is Hamiltonian for all
positive $k$ through $k=15$.  In this note we announce that $M_{33}$
and $M_{35}$ have Hamilton cycles.  The result was achieved by an
algorithmic improvement that made it possible to find a Hamilton path
in a reduced graph (of complementary necklace pairs) having 129,644,790
vertices, using a 64-bit personal computer.
\end{abstract}

\begin{keyword}
Hamilton cycles \sep middle levels \sep Boolean lattice \sep necklaces

\end{keyword}
\date{30 June 2006}
\end{frontmatter}

\section{Introduction}
\label{sec:intro}

\begin{figure}
\centering
\includegraphics[scale=0.7]{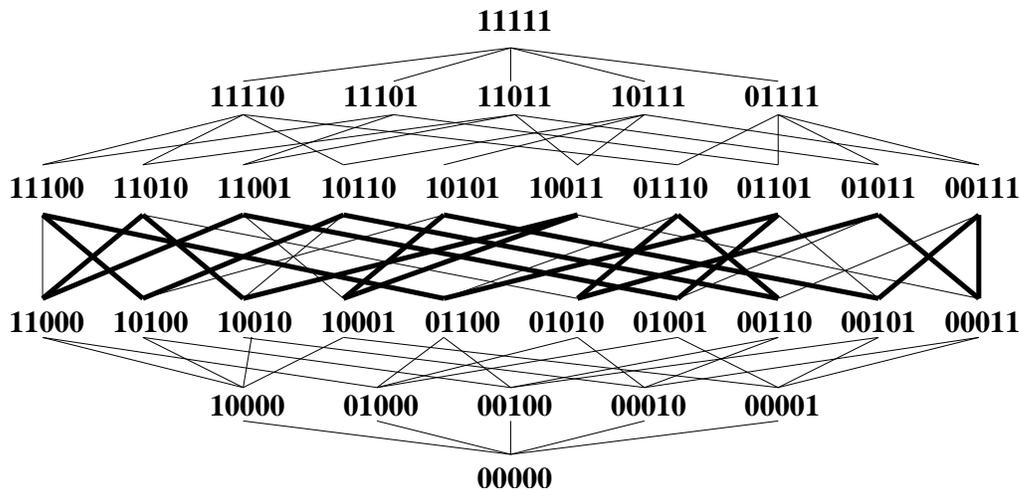}
\caption{ The Hasse diagram of $ {\mathcal B}_{5} $ showing a Hamilton
cycle in $M_5$. }
\label{mm2_5}
\end{figure}

Let $ {\mathcal B}_{n} $ be the $n$-atom Boolean lattice, i.e., the
partially ordered set of subsets of $[n]=\{1,2,\ldots,n\}$, ordered by
inclusion. The Hasse diagram of $ {\mathcal B}_{n} $ is isomorphic to
the $n$-cube, whose vertices are the $n$-bit binary numbers, with two
numbers adjacent if they differ in one bit position
(Figure~\ref{mm2_5}).  The $i$th level of ${\mathcal B}_n$ can thus be
viewed as the set of $n$-bit binary numbers with $i$ ones.

The \emph{\mlevels\ problem} is to determine if there is a Hamilton
cycle in the subgraph $M_{2k+1}$ of ${\mathcal B}_{2k+1}$ induced by
the \mlevels\ $k$ and $k+1$.  The heavier lines in Figure~\ref{mm2_5}
show a Hamilton cycle in $M_5$.  The graph $M_{2k+1}$ gained notoriety
as an example of a connected, vertex transitive graph, all of which
were conjectured by Lov\'{a}sz \cite{lovasz:1970:combinatorialst} to
have Hamilton paths.

The {\mlevels\ problem} remains open, in spite of the efforts of many
\cite{MR87a:05101,MR90d:05155,MR95b:05135,MR90a:05149,MR94m:05116,MR96h:05070}.
In 1990, in unpublished work, Moews and Reid verified that $M_{2k+1}$
is Hamiltonian for $1 \leq k \leq 11$.  In 1999, we extended this for
$12 \leq k \leq 15$ \cite{MR2000k:05179}.  In this note we announce
that $M_{33}$ and $M_{35}$ are Hamiltonian.  The new results are due
to an algorithmic improvement that made it possible to find a Hamilton
path in a reduced graph of complementary necklace pairs having
129,644,790 vertices, using a 64-bit personal computer.

It is known that any connected vertex transitive graph $G$ with $V(G)$
vertices has a cycle of length at least $\sqrt{3V(G)}$
\cite{MR80m:05059}.  For $M_{2k+1}$, the best absolute lower bound is
given by the following theorem in \cite{MR96c:94008}.
\begin{theorem}
\label{thm:savagewinkler}
The middle levels graph, $M_{2k+1}$, has a cycle of length at least
$$\left (1-\frac{\binom{2(j+1)}{j+1}}{2^{2(j+1)}}\right )V(M_{2k+1})$$
if for every $i \le j$, $M_{2i+1}$ has a Hamilton cycle.
\end{theorem}

Since our new results show $M_{2k+1}$ Hamiltonian for $1 \leq k \leq
17$, it follows that $M_{2k+1}$ has a cycle of length at least $0.867
V(M_{2k+1})$.

Recently, it was shown that $M_{2k+1}$ is ``asymptotically
Hamiltonian'' in the following sense \cite{MR2046083}: \emph{ There is a
constant $c$ such that for all $k$, $M_{2k+1}$ has a cycle of length
at least $(1-c/\sqrt{k}) V(M_{2k+1})$.}  In the other direction, it
has been shown that $M_{2k+1}$ has a closed spanning walk in which no
edge and no vertex occurs more than twice \cite{MR2195731}.

In Section~\ref{sec:middle2prob} we describe a reduction of the
problem which lessens the memory and computation requirements.  In
Section~\ref{sec:oldheur} we review the Hamilton cycle heuristic (SS)
from~\cite{MR2000k:05179}.  In Section~\ref{sec:newheur} we describe
the improved SSS heuristic which made the present results possible.

\section{Reducing the problem}
\label{sec:middle2prob}

We sketch here the reduction of the \mlevels\ problem that was was
used in the earlier work of Dejter~\cite{MR87a:05101} and of Moews and
Reid.

Given a string $x=x_1x_2 \ldots x_n$ of $n$ symbols, define the cyclic
shift $\sigma$ by $ \sigma(x_1x_2 \ldots x_n)=x_2x_3 \ldots x_nx_1.$
Let $\sigma^{1}(x)=\sigma(x)$ and for $i \geq 0$, $\sigma^{i+1}(x)=
\sigma(\sigma^{i}(x))$.  Define a relation $ \sim $ on the set of
$n$-bit binary numbers (regarded as $n$-bit strings) by $x \sim y $ if
and only if $y=\sigma^i(x)$ for some integer $i$.  The relation $ \sim
$ is an equivalence relation and the equivalence classes are called
\emph{necklaces}.  Denote the necklace of $x$ by $\nu(x)$.

Fix $n=2k+1$.  We use necklaces to define a quotient graph of the
middle levels graph $M_n$.  Let $N_n$ be the graph whose vertices are
the necklaces of the vertices of $M_n$ (i.e., necklaces of $2k+1$-bit
strings with $k$ or $k+1$ ones).  The edges of $N_n$ are those pairs
$\nu(x)\nu(y)$ such that $xz \in E(M_n)$ for some $z \in \nu(y)$.
Note that the necklace of a $2k+1$-bit binary number with $k$ or $k+1$
ones has exactly $2k+1$ elements, so $N_n$ is smaller than $M_n$ by a
factor of $n$.

The \emph{complement} $\overline{b}$ of a binary digit $b$ is 1 if
$b=0$ and 0 if $b=1$.  Extend this to binary strings by bitwise
complement and to necklaces by by $\overline{\nu{(x)}} =
\nu{(\overline{x})}$.  Note that this is well-defined since $y =
\sigma^i(x)$ if and only if $\overline{y} = \sigma^i(\overline{x})$.

We use these complementary necklace pairs to further reduce the
problem size by a factor of 2 by observing that $\nu(x)\nu(y) \in
E(N_n)$ if and only if $\overline{\nu(x)}\,\overline{\nu(y)} \in
E(N_n)$.  Now define an equivalence relation
$\stackrel{\diamond}{\sim}$ on the necklaces in $V(N_n)$ by $X
\stackrel{\diamond}{\sim} Y$ if either $X = Y$ or $X = \overline{Y}$
and denote the equivalence class of $\stackrel{\diamond}{\sim}$
containing $X$ by $\rho(X)$.  Since every string in $X$ has odd
length, $X \not = \overline{X}$, and every equivalence class $\rho(X)$
has exactly 2 elements.  Construct the reduced graph, $R_n$, whose
vertices are the equivalence classes $\{\rho(X) \ | \ X \in V(N_n) \}$
with edges $\rho(X)\rho(Y) \in E(R_n)$ if $XY \in E(N_n)$ or
$X\overline{Y} \in E(N_n)$.  Observe that if $Z=\nu(0^k1^{k+1})$ then
$\overline{Z} = \overline{\nu(0^k1^{k+1})} = \nu(1^k0^{k+1})$, and so
$Z\overline{Z} \in E(N_n)$. Hence $\rho(Z)\rho(\overline{Z}) \in
E(R_n)$ and $R_n$ has loops, so it is not a simple graph.

\begin{figure}
\centering
\subfigure[$R_5$ path.]{%
\includegraphics[scale=1.0]{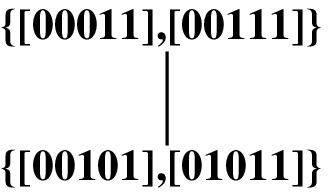}}\quad
\subfigure[Using loops to lift complementary paths to cycle in $N_5$.]{%
\includegraphics[scale=1.0]{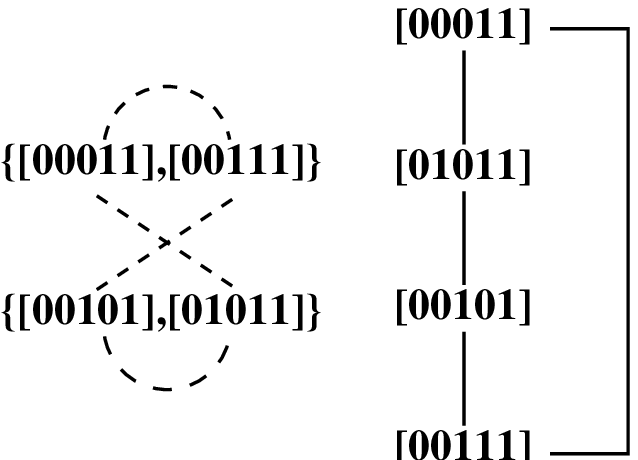}}\quad
\subfigure[Lifting cycles in $N_5$ to cycle in $M_5$.]{%
\includegraphics[scale=1.0]{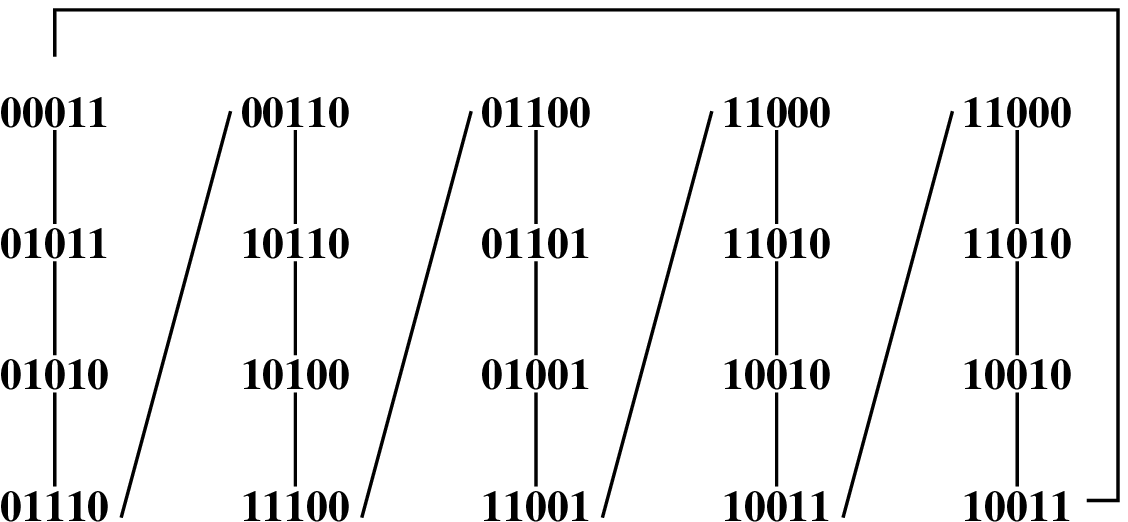}}
\caption{ Lifting a path in $R_5$ to a cycle in $N_5$ and then a cycle in $M_5$. }
\label{reductionfig}
\end{figure}

We exploit the fact that $R_n$ has loops to show that a Hamilton path
from the distinguished vertex $r_1= \rho(\nu(0^k1^{k+1}))$ to the
distinguished vertex $r_l=\rho(\nu(0(01)^k))$ in $R_n$ can be used to
construct a cycle in $N_n$ which can be lifted to a cycle in
$M_n$. The path in $R_n$ gives rise naturally to a pair of paths in
$N_n$ where corresponding pairs of vertices from each path are
complements of each other. Because the distinguished vertices in $R_n$
each have incident loops, these paths link to form a cycle in
$N_n$. Finally, since $k$ and $2k+1$ are relatively prime, we can use
suitably chosen necklace representatives to extend the cycle in $N_n$
to a cycle in $M_n$.  The case for $M_5$ is illustrated in
Figure~\ref{reductionfig}, where for $x \in V(M_5)$, $[x]$ denotes the
necklace $\nu(x) \in V(N_5)$ and for necklace $X \in V(N_5)$, $\{X,
\overline{X}\}$ denotes $\rho(X)$ in $V(R_5)$.

\section{The Hamilton cycle heuristic}
\label{sec:oldheur}

Given a graph $G$ and vertices $s,t \in V(G)$ we would like to find,
if possible, a Hamilton path starting at $s$ and ending at $t$. 

\begin{figure}
\centering\includegraphics[width=\linewidth]{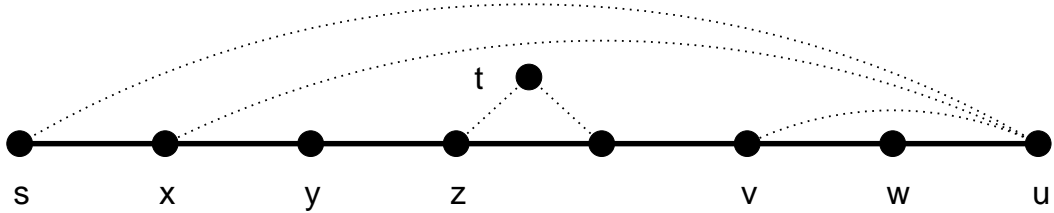}
\caption{Backtrack search from $s$ has reached $u$ and cannot
continue.}
\label{rot1}
\end{figure}

A standard backtrack search attempts to construct such a path by
starting at $s$ and extending the path to a new vertex as long as
possible.
Whenever it is impossible to further extend the path from a vertex
$x$, the search ``backs up'' to the predecessor, $y$, of $x$ on the
path and the path is extended (if possible) from $y$ to one of its
other neighbors.
Figure~\ref{rot1}
illustrates a path $P$ from a starting vertex $s$ to a vertex $u$
whose only neighbors in $G$ are $s$, $x$ and $v$, which are already on
$P$.  In this case, backtrack search would back up to $w$, then $v$
and eventually back to $z$, at which time $t$ would be added to the
path.

\begin{figure}
\centering
\subfigure[Path ending at $y$ after rotation at $x$.]{%
\label{bad1:a}%
\includegraphics[width=0.45\textwidth]{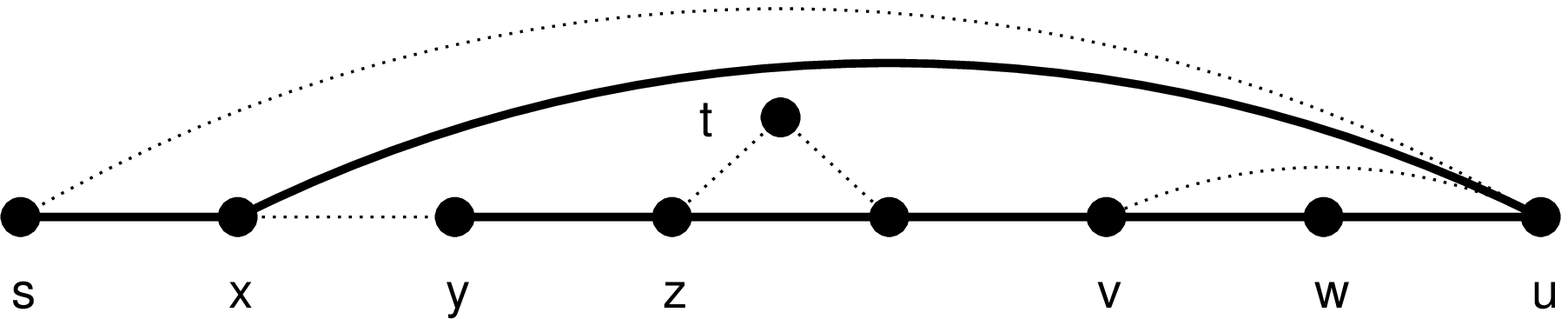}}\quad
\subfigure[Path ending at $x$ after rotation at $s$.]{%
\label{bad1:b}%
\includegraphics[width=0.45\textwidth]{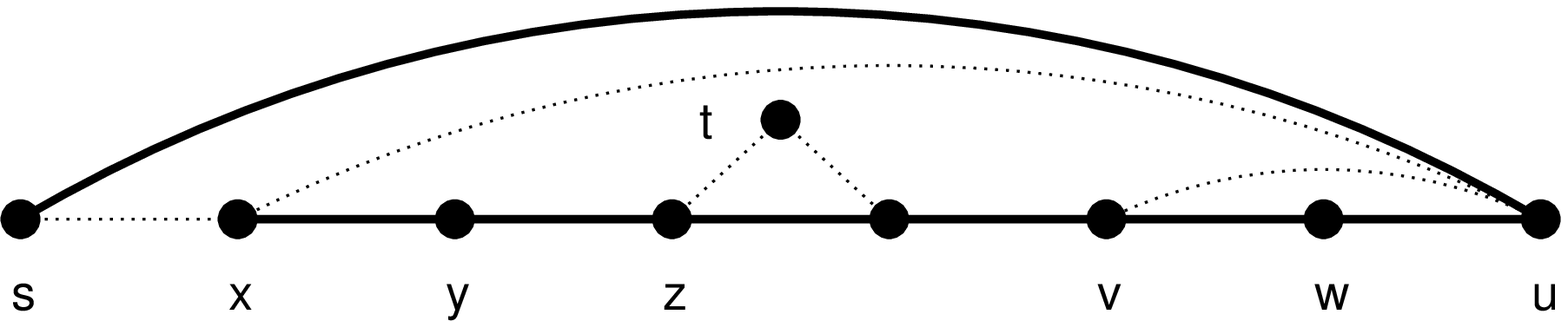}}\quad
\subfigure[Path ending at $w$ after rotation at $v$.]{%
\label{bad1:c}%
\includegraphics[width=0.45\textwidth]{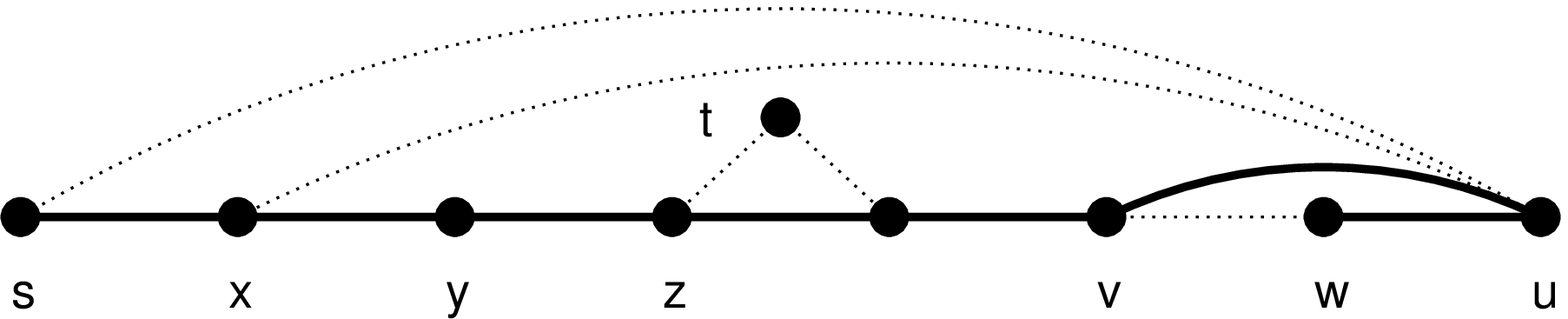}}\quad
\subfigure[Path ending at $w$ after rotation at $s$ followed by rotation at $u$.]{%
\label{bad1:d}%
\includegraphics[width=0.45\textwidth]{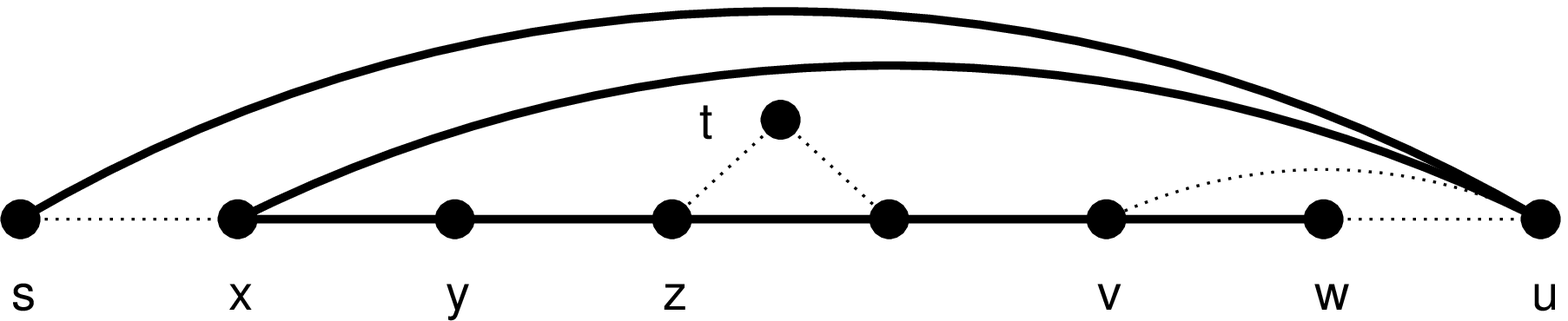}}
\caption{ Possible endpoints after one or more rotations }
\label{bad1}
\end{figure}

\Posa~\cite{MR52:10497} observed that reaching a dead end in backtrack
search could be used as an opportunity to modify the current path by a
\emph{rotation} and thus possibly to continue. If
$P=(u_1,u_2,\ldots,u_k)$ is a path in $G$ and there is an edge
$u_ku_j$ for some $0<j<k$, then the \emph{rotation of $p$ at $u_j$} is
the path $P'=(u_1,u_2,\ldots u_j,u_k, u_{k-1},\ldots, u_{j+1})$
obtained by deleting edge $u_ju_{j+1}$ from $P$ and inserting
$u_ju_k$.  Path $P'$ has the same length as $P$, but a different
endpoint.  Figure~\ref{bad1:a} shows the rotation at $x$ of the path
in Figure~\ref{rot1}.

Define \emph{ PosaSearch} to be the Hamilton path heuristic which
constructs a path $P$ by starting at a vertex $s$ and iterating the
following: extend $P$, avoiding $t$ until the end, until no longer
possible; When $P$ cannot be extended from an endpoint $u$, select a
neighbor $v$ of $u$ and perform a rotation of $P$ at $v$.  PosaSearch
may not succeed in finding a Hamilton path even if one exists.
Furthermore, it may run indefinitely.

For the path $P$ in Figure~\ref{rot1}, PosaSearch can only transform
$P$ into one of the paths shown in Figure~\ref{bad1}.  None of these
transformations will ever allow the vertex, $t$, to be added to the
path.

Rotations transform a path and alter the order of vertices on it.  For
a path $P$, let $P(i)$ be the $i$th vertex of $P$, let
$\mathit{pos}(P,v)$ be the position of vertex $v$ on $P$, and let
$|P|$ be the number of vertices on $P$.  If $P$ is a path with
endpoint $u$ and if $P(j)$ is adjacent to $u$, we can rotate $P$ at
$P(j)$ to arrive at a new path $P'$.  The position of a vertex $v$ on
$P'$ is then given by
\begin{equation}
\mathit{pos}(P',v) = \left\{ \begin{array}{ll}
\mathit{pos}(P,v) & \mbox{ if $\mathit{pos}(P,v) \leq j$} \\
|P|-\mathit{pos}(P,v)+j+1 & \mbox{ otherwise. }
\end{array}
\right.
\label{eq:newpos}
\end{equation}
Similarly, the vertex now in position $i$ on $P'$ is
\begin{equation}
P'(i) = \left\{ \begin{array}{ll}
P(i) & \mbox{ if $i \leq j$} \\
P( |P|-i+j+1) & \mbox{ otherwise. }
\end{array}
\right.
\label{eq:newsucc}
\end{equation}

The SS heuristic from \cite{MR2000k:05179} uses a variation of
PosaSearch that looks ahead before performing any rotation.  SS
extends the path as is done in PosaSearch until all neighbors of the
endpoint are already on the path.  If the path is not already a
Hamilton path, it then uses breadth-first search (BFS), and repeated
application of eqs.  (\ref{eq:newpos}) and (\ref{eq:newsucc}) to
search for a sequence of one or more rotations guaranteed to result in
a path which can be further extended.  If such a sequence is found,
the sequence of rotations is performed and the path is extended.  If
no such sequence exists, the SS heuristic terminates without having
found a Hamilton path.  Figure~\ref{bad1} shows four paths obtainable
by one or more rotations from the path shown in Figure~\ref{rot1}, In
this example, there are only four possible endpoints, $u$, $v$, $w$,
and $x$. The path cannot be extended from any of these, so SS would
terminate after evaluating all four.

\section{The new SSS heuristic}
\label{sec:newheur}

The original work on the SS heuristic \cite{MR2000k:05179} used a
400MHz Intel Pentium-II system with 192MB RAM. Results were later run
on a 2.4GHz Intel Pentium 4 system with 512MB of RAM
(See~\cite{Shields:2004:PhD}).We recently acquired an AMD Athlon 3500+
system with 2GB of RAM and converted the program to use 64-bit values
for most internal computations. With this system and the SS heuristic
we were able to find a Hamilton path in $R_{33}$ in about 3.5
weeks. The program ran entirely in memory, using approximately 1GB of
RAM.  However, it was unlikely that $R_{35}$ would be feasible.

Therefore, using performance profiling tools, we analyzed the code
performance and discovered that, contrary to expectations, a large
part of the time was spent on performing the rotation operations,
rather than on the BFS to find promising sequences of rotations.  As a
result, we made the following changes to the heuristic that resulted
in dramatic speed improvements.

First, when BFS finds a sequence of rotations that will enable an
extension of the path, instead of actually performing each rotation of
the path the SSS heuristic encodes the sequence of rotations as a list
of ordered pairs, each representing the number of vertices on the path
at time of rotation and the rotation point.  This list, together with
repeated application of eqs.  (\ref{eq:newpos}) and
(\ref{eq:newsucc}), suffices to calculate the actual vertex in a given
position or the current position of a given vertex, when the sequence
of rotations represented by the list of ordered pairs has not been
explicitly performed.

Secondly, the SSS heuristic periodically (but rarely) goes ahead and
performs all the rotations encoded by the list of ordered pairs.
Given a path $v_0,\ldots,v_l$, if we perform a rotation at $v_i$ to
make $v_{i+1}$ the new endpoint, then there is a section of $i+1$
vertices up to $v_i$ that will not be moved and a section of $l-i$
vertices after $v_i$ that will be reversed. We may describe each of
these sections (or \emph{blocks}) of the path by a triple,
$\{s,n,d\}$, where $s$ is the index of the first vertex in the block,
$n$ is the number of vertices in the block and $d$ is a direction
indicating whether the block is now in original (forward) order or
reversed. Our single rotation at $v_i$ is represented by the triples
$\{0,i+1,\rightarrow\}$ and $\{i+1, l-i,\leftarrow \}$. As we process
a rotation from the list of ordered pairs we first add a block in the
forward direction, representing any vertices added by extension since
the last rotation. Next we determine the block containing the new
endpoint and (possibly) split it into two blocks. All blocks before
the new endpoint remain unchanged, while the sequence of blocks after
the new endpoint, as well as the direction indicator of each such
block, is reversed. After the list of blocks has been created
representing all the stored rotations, the path can be copied block by
block to effect the series of rotations in a single copy operation.

These modifications allow a sequence of rotations to be accumulated,
without being performed, at the relatively smaller cost of increasing
time for calculation of vertex position or vertex in a given position.
The accumulated sequence of rotations is performed on a schedule
discussed in the next subsection.  This significantly reduced the
average work per rotation.

Figure~\ref{blockfig} illustrates the process for a list of three
rotations, as seen in one run of the program. The subscripts in each
block indicate the order of creation of the block. For computational
convenience, a rotation point, $v$, is actually represented in the
list of rotations by its position on the path.

\begin{figure}
{\small
\begin{tabular}{ll}
\multicolumn{2}{l}{\textbf{Saved list of three rotations: $\{35,15\}$, $\{36,9\}$, $\{40,32\}$}}\\
\multicolumn{2}{l}{1.  $\{35,15\}$ - 35 vertices on path, rotation at $v_{15}$}\\
& Create initial block 0\\
&  \fbox{0,35,$\rightarrow_0$} \\ 
& Split block 0, into blocks 0 and 1\\
& \fbox{0,16,$\rightarrow_0$} \fbox{16,19,$\rightarrow_1$}\\ 
& Rotate block 1\\
& \fbox{0,16,$\rightarrow_0$} \fbox{16,19,$\leftarrow_1$}\\ 
\multicolumn{2}{l}{2.  $\{36,9\}$ - 36 vertices on path, rotation at $v_9$}\\
& Add block 2\\
& \fbox{0,16,$\rightarrow_0$} \fbox{16,19,$\leftarrow_1$} \fbox{35,1,$\rightarrow_2$}\\ 
& Split block 0, into blocks 0 and 3\\
& \fbox{0,10,$\rightarrow_0$} \fbox{10,6,$\rightarrow_3$} \fbox{16,19,$\leftarrow_1$} \fbox{35,1,$\rightarrow_2$}\\ 
& Rotate blocks 3 through 2\\
& \fbox{0,10,$\rightarrow_0$} \fbox{35,1,$\leftarrow_2$} \fbox{16,19,$\rightarrow_1$}  \fbox{10,6,$\leftarrow_3$} \\ 
\multicolumn{2}{l}{3.  $\{40,32\}$ - 40 vertices on path, rotation at $v_{32}$}\\
& Add block 4\\
& \fbox{0,10,$\rightarrow_0$} \fbox{35,1,$\leftarrow_2$} \fbox{16,19,$\rightarrow_1$}  \fbox{10,6,$\leftarrow_3$}  \fbox{36,4,$\rightarrow_4$} \\ 
& Split block 3, into blocks 3 and 5\\
& \fbox{0,10,$\rightarrow_0$} \fbox{35,1,$\leftarrow_2$} \fbox{16,19,$\rightarrow_1$}  \fbox{13,3,$\leftarrow_3$}  \fbox{10,3,$\leftarrow_5$}  \fbox{36,4,$\rightarrow_4$} \\ 
& Rotate blocks 5 and 4 \\
& \fbox{0,10,$\rightarrow_0$} \fbox{35,1,$\leftarrow_2$} \fbox{16,19,$\rightarrow_1$}  \fbox{13,3,$\leftarrow_3$}  \fbox{36,4,$\leftarrow_4$}  \fbox{10,3,$\rightarrow_5$} \\ 
\multicolumn{2}{l}{\textbf{ Final path:}}\\
& $\underbrace{v_0,\ldots,v_9}_0$,$\underbrace{v_{35}}_2$,$\underbrace{v_{16},\ldots,v_{34}}_1$,$\underbrace{v_{15},v_{14},v_{13}}_3$,$\underbrace{v_{39},v_{38},v_{37},v_{36}}_4$,$\underbrace{v_{10},v_{11},v_{12}}_5$
\end{tabular}
}

\centering
\caption{ Performing a sequence of rotations from a list. }
\label{blockfig}
\end{figure}

\subsection{Graph representation and storage}
\label{sec:m-graphrep}

In the middle levels of ${\mathcal B}_{35}$, the reduced graph,
$R_{35}$, has over 129 million vertices. An array of 32-bit integers
with one entry per vertex takes approximately half a gigabyte of
memory. Our system had only 2GB of RAM, so storage for basic
information on vertices and paths was limited. In this environment, it
is infeasible to store adjacency lists, so we recalculate adjacency
lists as needed rather than storing them
(see~\cite{Shields:2004:PhD}). We reduced all arrays to the smallest
native size (char, short, int, or long) that was sufficient for the
data and still experienced unacceptable levels of memory swapping.

The vertices of $R_{35}$ are represented by the lexicographically
least elements (as sets) of each $35$-bit necklace with $17$ ones.
Thus the two higher order bits will always be 0 and we need only use
33 bits for the internal representation. However, we store the lower
32 bits in 32-bit integer arrays and add the high-order bit according
to the position of a vertex in the array. This saves half a gigabyte
of RAM as compared to using 64-bit long integers. We still use 64-bit
long integers for internal calculations, but none of our major storage
arrays needs larger than 32-bit entries.

When copying the path to execute a series of rotations, we copy it
to the position array (Pos) and then rebuild the position array after
the copy as this results in less swapping than using memory that has
not recently been referenced, such as the parent array used in
building a BFS tree.

Fine tuning in this way resulted in the program using approximately
85\% to 90\% of the CPU for most of the run while finding a path in
$R_{35}$. Somewhat higher CPU usage was observed toward the end of the
run, even though BFS trees, with their storage requirement, were being
computed much more frequently. 

\subsection{Performance of the heuristic on the \mlevels\ graph}
\label{sec:performance}

The SSS heuristic was applied to search the reduced middle levels
graphs $R_{2k+1}$, $k = 1,2, \ldots 17$, to try to find a Hamilton
path from vertex $\rho(\nu(0^{k+1}1^k))$ to vertex
$\rho(\nu(0(01)^k))$, and therefore a Hamilton cycle in $M_{2k+1} $.
A Hamilton path meeting these requirements was found for each $k \leq
17$.

The results are summarized in Table~\ref{runtimes}.  We measured
elapsed time using a timer with a resolution of 1 second. Runs that
start and complete in the same second show a time of 0
seconds. Table~\ref{runtimes} also shows earlier results obtained with
the SS heuristic on a 2.4GHz Intel Pentium 4 system with 512MB of RAM.
Many factors affect the difference in performance between the two
systems, including operating system and compiler differences as well
as processor, RAM, and disk speed. The newer system ran the SS
heuristic approximately twice as fast as the older system. Thus for
$k=15$, the largest value for which we have both results, we see that
hardware doubled the speed and the switch to the SSS algorithm
achieved a further 16-fold increase.

\begin{table}[ht]
\centering
\caption{Running time to find a Hamilton cycle in the middle 
levels graph}
\label{runtimes}

\medskip

\begin{tabular}{|c|c|c|c|c|c|} \hline

$k$   & $n=2k+1$  & \# vertices in $R_n$ & 
\# vertices in $M_n$ & \multicolumn{2}{c|}{Time (Secs)}  \\ 
    &    &           &              &  \multicolumn{1}{c|}{SS} & SSS\\ 
    &    &           &              &  \multicolumn{1}{c|}{(32-bit)} & (64-bit) \\ \hline
8   & 17  & 1,430  &48,620  & 0 & 0\\
9   & 19 & 4,862 & 184,756  & 0 & 0\\
10  & 21 & 16,796  & 705,432  & 1 & 0\\
11  & 23 & 58,786 & 2704,156  & 5 & 2\\
12  & 25 & 208,012 & 10,400,600 & 105 & 10 \\
13  & 27 & 742,900 & 40,116,600  & 1,732 & 99 \\
14  & 29 & 2,674,440 & 155,117,520  & 24,138 & 799\\
15  & 31 & 9,694,845 & 601,080,390  & 307,976 & 9,446\\ 
    &    &           &              &  (3.6 days) & \\ 
16  & 33 & 35,357,670 &  2,333,606,220  & - & 106,118 \\ 
    &    &           &              & & (1.2 days)\\ 
17  & 35 & 129,644,790  & 9,075,135,300  & - & 1,765,497 \\ 
    &    &           &              & &  (20.4 days)\\ \hline

\end{tabular}
\end{table}

As noted earlier, the deferral of rotations also involves a cost. We
experimented with the number of deferred rotations using values
ranging from $ \log V(G)$ to $ \sqrt{V(G)}$.  Our results were
obtained using a value of $ \sqrt{V(G)}$, or 11,387 for the reduced
graph, $R_{35}$.  There is some evidence to suggest that a small
performance improvement could be obtained by varying this number
during the running of the program, particularly in very large graphs.

\bibliographystyle{elsart-num}






\end{document}